\documentclass[12pt,a4paper,twoside]{article}
\usepackage{indentfirst}
\usepackage{amsmath}
\usepackage{amssymb}
\usepackage{t1enc}
\usepackage[latin2]{inputenc}
\usepackage[english]{babel}
\usepackage[amsthm]{ntheorem}
\usepackage{graphicx}
\usepackage{psfrag}
\usepackage{longtable}
\usepackage{picinpar}

\renewenvironment{proof}{\noindent\textit{Proof}.}{\hfill $\square$}
\newenvironment{proof2}{\noindent\textit{Proof of Theorem \ref{boxmain}}.}{\hfill $\square$}
\newenvironment{proof1}{\noindent\textit{Proof of Theorem \ref{main}}.}{\hfill $\square$}
\newenvironment{proofx}{\noindent\textit{Proof}.}{}

\newcommand*{\R}{\ensuremath{\mathbb R}}
\newcommand*{\C}{\ensuremath{\mathbb C}}
\newcommand*{\Z}{\ensuremath{\mathbb Z}}
\renewcommand*{\P}{\ensuremath{\mathbb P}}

\newcommand*{\sd}{\ensuremath{\textup{sd}}}
\newcommand*{\N}{\ensuremath{\textup{N}}}

\newcommand*{\id}{\ensuremath{\textup{id}}}

\newcommand*{\ind}{\ensuremath{\textup{ind}}}
\newcommand*{\dni}{\ensuremath{\textup{coind}}}
\newcommand*{\B}{\ensuremath\textnormal{B}}

\newcommand*{\susp}{\ensuremath{\textup{susp}}}
\newcommand*{\cone}{\ensuremath{\mathcal{C}}}
\newcommand*{\CN}{\ensuremath{\textup{CN}}}
\newcommand*{\con}{\ensuremath{\textup{connectivity}}}

\renewcommand*{\H}{\ensuremath{\mathcal H}}

\newcommand{\subpart}[1]{\noindent\textbf{#1:}}

\theoremseparator{.}
\theoremstyle{plain}
\newtheorem{thm}{Theorem}
\newtheorem{lemma}[thm]{Lemma}

\newtheorem{def2}[thm]{Definition}

\newtheorem{constr}[thm]{Construction}

\theoremstyle{remark}\newtheorem{rem}[thm]{Remark}
\theoremstyle{definition}\newtheorem{Def}[thm]{Definition}

\newlength{\Uj}
\newcommand{\chr}[1]{{\raisebox{\Uj}{$\chi$}}{(#1)}}

\begin{document}
\title{Homotopy types of box complexes}
\author{P\'eter Csorba\footnote{\ Supported
        by the joint Berlin/Z\"{u}rich graduate
        program Combinatorics, Geometry, and Computation,
        financed by ETH Z\"{u}rich and the German Science
        Foundation (DFG).}\\
{\footnotesize  Institute of Theoretical Computer Science}\\[-1.5mm]
{\footnotesize  ETH Zürich, CH-8092 Zürich, Switzerland}\\[-1.5mm]
{\footnotesize E-mail: \texttt{pcsorba@inf.ethz.ch}}}
\maketitle

\begin{abstract}
In \cite{MZ} Matou\v{s}ek and Ziegler compared various
topological lower bounds for the chromatic number.
They proved that Lov\'{a}sz's original bound \cite{L}
can be restated as $\chr G \geq \ind (\B(G)) +2$.
Sarkaria's bound \cite{S} can be formulated as
$\chr G \geq \ind (\B_0(G)) +1$. It is known that these lower
bounds are close to each other, namely the difference between them
is at most 1. In this paper we study these lower bounds, and the homotopy
types of box complexes. Some of the results was announced in \cite{MZ}.
\end{abstract}

\section{Introduction}
In \cite{MZ} Matou\v{s}ek and Ziegler compared various
topological lower bound for the chromatic number.
They reformulated Lov\'{a}sz's original bound \cite{L} and
Sarkaria's bound \cite{S} in terms of various box complexes:

\begin{thm}[The Lov\'asz bound \cite{MZ}]\label{lov}
For any graph $G$
\[\chr G \geq \ind (\B(G)) +2.\]
\end{thm}
\begin{thm}[The Sarkaria bound \cite{MZ}]\label{sark}
For any graph $G$
\[\chr G \geq \ind (\B_0(G)) +1.\]
\end{thm}

\noindent
We will study these lower bounds in this paper, which is organized as follows.

Section \ref{sec2} contains the definition
of the box complexes of graphs and we fix some notation.

In Section \ref{sec3} we prove that the box complex $\B_0(G)$ is
$\Z_2$-homotopy equivalent to the suspension of $\B(G)$. This makes
the connection between these two bounds explicit.
Since $\ind(X)\leq\ind(\susp(X))\leq\ind(X)+1$ the difference between
the right side of Lov\'asz and the Sarkaria bound is at most one.

From topological point of view it is possible that these two bounds
are not the same. We construct a $\Z_2$-space $X_{2h}$ such that
$\ind(\susp(X_{2h}))=\ind(X_{2h})$ in Section \ref{sec6}.

However we need a graph such that its box complex $\B(G)$ has this property.
In Section \ref{sec4} we show that the homotopy type of the box complex
$\B(G)$ (which is homotopy equivalent to the neighborhood complex)
can be "arbitrary"; in Section \ref{sec5} we extend this result to
$\Z_2$-homotopy equivalence. This allows us to construct a graph $G$ such
that the gap between these two bounds is 1.
This means that the Lov\'asz bound can be better than the Sarkaria bound,
which answers a question of Matou\v{s}ek and Ziegler \cite{MZ}.

Finally in Section \ref{sec7} we show that both of these topological lower
bounds can be arbitrarily bad. Our examples are purely topological.


\section{Preliminaries}\label{sec2}

In this section we recall some basic facts of graphs and simplicial
complexes and topology to fix notation. The interested reader
is referred to~\cite{M} or~\cite{B} and~\cite{Hat01} for details.

\subpart{Graphs}
Any graph $G$ considered will be assumed to be finite, simple, connected,
and undirected, i.e. $G$ is given by a finite set $V(G)$ of
\emph{vertices} and a set of \emph{edges} $E(G) \subseteq \binom{V(G)}{2}$.
A \emph{graph coloring with $n$ colors} is a homomorphism
$c: G\to K_n$, where $K_n$ is the complete graph on $n$ vertices and
the \emph{chromatic number} $\chr G$ of $G$ is the
smallest $n$ such that there exists a graph coloring of $G$ with $n$ colors.
The {\em common neighbor} of $A\subseteq V(G)$ is
$\CN(A)\ =\ \{v\in V(G)\colon \{a,v\}\in E(G)\mbox{ for all } a\in A\}$.
For two disjoint sets of vertices $A,B \subseteq V(G)$ we define
$G[A,B]$ as the (not necessarily induced)
subgraph of $G$ with $V(G[A,B])= A \cup B$ and
$E(G[A,B])=\{ (a,b) \in E(G) \colon a\in A,b\in B\}$.

\subpart{Simplicial Complexes}
A \emph{simplicial complex} $\mathcal{K}$ is a finite hereditary set system.
We denote its vertex set by $V(\mathcal{K})$ and its barycentric
subdivision by $\sd(\mathcal{K})$.

\noindent
For sets $A,B$ define
$ A\uplus B := \{(a,0)\colon a\in A\} \cup \{(b,1)\colon b\in B\}$.

\subpart{Neighborhood Complex}
The \emph{neighborhood complex} is
$\N(G)= \{S\subseteq V(G)\colon \CN(S)\neq\emptyset\}$

\subpart{Box Complex}
The \emph{box complex} $\B(G)$ of a graph $G$
(the one introduced by Matou{\v s}ek and Ziegler \cite{MZ}) is defined by
\[
\B(G)  := \{ A\uplus B \colon A,B\subseteq V(G),\ 
                         A\cap B = \emptyset,\ 
                         G[A,B] \text{ is complete bipartite},\ 
                         \CN(A)\neq\emptyset\neq\CN(B)
                       \}
\]
The vertices of the box complex are $V_{1} :=\left\{v\uplus\emptyset
\colon v\in V(G)\right\}$ and $V_{2} :=\left\{\emptyset\uplus v\colon
v\in V(G)\right\}$ for all vertices of $G$. The
subcomplexes of $\B(G)$ induced by $V_{1}$ and $V_{2}$ are disjoint
subcomplexes of $\B(G)$ that are both isomorphic to the neighborhood
complex $\N(G)$. We refer to these two copies as \emph{shores} of
the box complex. The box complex is endowed with a $\Z_{2}$-action
which interchanges the shores.


A \emph{different} box complex $\B_0 (G)$:
\[
\B_0 (G)\ = \left\{ A\uplus B \colon
A,B\subseteq V(G),\ A\cap B=\emptyset,\ G[A,B]\mbox{ is
  complete bipartite}\right\}
\]
The \emph{cones over the sores} complex $\B_{\cone}(G)$
(only for technical reason):
\[
\B_{\cone} (G){=}\B(G)\cup \left\{(x, A\uplus\emptyset)
\colon A\subseteq V(G), \CN(A){\neq}\emptyset\right\}
\cup \left\{(\emptyset\uplus B, y)
\colon B\subseteq V(G), \CN(B){\neq}\emptyset\right\},
\]
where we assume that $x,y\not\in V(G)$.
($\B(G),\B_0(G),\B_{\cone} (G)$ are $\Z_2$-spaces.)


\medskip
\subpart{$\Z_2$-space}
A \emph{$\Z_2$-space} is a pair $(X,\nu)$ where $X$ is a topological
space and $\nu\colon X\to X$, called the $\Z_2$-action,
is a homeomorphism such that
$\nu^2=\nu\circ\nu=\mbox{\rm{id}}_X$. If $(X_1,\nu_1)$ and $(X_2,\nu_2)$
are $\Z_2$-spaces, a \emph{$\Z_2$-map} between them is
a continuous mapping $f\colon X_1\to X_2$ such that $f\circ\nu_1=\nu_2\circ f$.
The sphere $S^n$ is understood as a $\Z_2$-space with the
antipodal homeomorphism $x\to -x$.
We will consider only finite dimensional \emph{free} $\Z_2$-complexes
(free means that the $\Z_2$-action $\nu$ has no fixed point).

\subpart{$\Z_2$-index}
We define the \emph{$\Z_2$-index} of a $\Z_2$-space $(X,\nu)$ by
\[
\ind(X)= \min \left\{n\geq 0 \colon \mbox{ there is
a $\Z_2$-map $(X,\nu)\to (S^n,-)$}\right\}
\]
(the $\Z_2$-action $\nu$ will be omitted from the notation if it is clear
from the context).
The Borsuk--Ulam Theorem can be re-stated as \  $\ind(S^n)=n$.

Another index-like quantity of a $\Z_2$-space, the \emph{dual index}
can be defined by
\[
\dni(X)= \max \left\{ n\geq 0 \colon \mbox{ there is
a $\Z_2$-map $S^n \stackrel{\Z_2}{\longrightarrow} X$}\right\}.
\]

The consequence of the
Borsuk--Ulam Theorem is that $\dni(X)\leq\ind(X)$. We call a free
$\Z_2$-space \emph{tidy} if $\dni(X)=\ind(X)$.
(In general $\ind(X)\geq\dni(X) \geq \con(X)+1$ \cite{M}.)

\medskip
A $\Z_2$-map $f\colon X\to Y$ is a $\Z_2$-equivalence if
there exist a $g\colon Y\to X$ such that $g\circ f$ and
$f\circ g$ are homotopic to $\id_X$ and $\id_Y$ respectively.
A general reference for $\Z_2$-spaces is \cite{Bre67}.

\section{The connection between $\B_{\cone}(G)$, $\B_0(G)$ and $\B(G)$}\label{sec3}

In this section we will prove that $\B_0(G)$ and $\susp(\B(G))$ are
$\Z_2$-homotopy equivalent. The reason is that the box complex is 'nearly'
$\N(G)\times [0,1]$.

\begin{rem}\label{r}
One can use Lov\'{a}sz's bound 
to prove Kneser's Conjecture \cite{K}.
The box complexes of Kneser graphs (Schrijver graphs) are tidy
spaces \cite{L} (
spheres up to homotopy \cite{BL}).
This means that one can prove Kneser's Conjecture by using Sarkaria's
bound 
(or any higher suspension) as well.
\end{rem}

\begin{lemma}
$\B_{\cone}(G)$ is $\Z_2$-homotopy equivalent to $\B_0(G)$.
\end{lemma}
\begin{proof}
$\B_{\cone}(G)$ was obtained from $\B(G)$ by attaching two cones $C_1,C_2$
over the shores, while $\B_0(G)$ is $\B(G)$ plus
two simplices $\Delta_1,\Delta_2$ covering the shores.

We consider the following two quotient CW-complexes.
$(\B_{\cone}(G)/C_1)/C_2$ and
$(\B_0(G)/\Delta_1)/\Delta_2$ (the order of the factorization does not
matter since we collapse disjoint subcomplexes).
It is obvious that they are the same CW-complexes and since $C_i,\Delta_i$
are contractible spaces $\B_{\cone}(G)$ and $\B_0(G)$ are $\Z_2$-homotopy
equivalent.
\end{proof}
\begin{lemma}\label{lemma2}
$\B_{\cone}(G)$ is $\Z_2$-homotopy equivalent to $\susp(\B(G))$.
\end{lemma}
\begin{proofx}
$\B_{\cone}(G)$ is a subcomplex of $\susp(\B(G))$. The idea of the proof
is to start with $\susp(\B(G))$, and get rid of the extra simplexes one by
one (using deformation retraction) such that finally we get $\B_{\cone}(G)$.
We will work with one cone (half) of the suspension. Since we want a
$\Z_2$-retraction, on the other cone we have to do the $\Z_2$-pair of
each step.

Let $x$ be the apex of the cone over the first shore in $\susp(\B(G))$
($y$ is the other apex).
We will define (by induction) sequences of simplicial complexes such that
\[
\susp(\B(G))=:X_0\supset X_1\supset\dots \supset X_N=\B_{\cone}(G)
\]
and $X_{i+1}$ is a $\Z_2$-deformation retraction of $X_i$.

Let assume that we already defined $X_n$. We choose a simplex
$\sigma\in X_n$ such that
\begin{enumerate}
\item $x\in\sigma$, and the rest of the vertices of $\sigma$ are
      from the second shore,
\item no other simplex in $X_n$ containing $x$ has \emph{more}
      vertex from the second shore,
      and it has at least one vertex from the second shore.
\end{enumerate}
The vertex set of $\sigma$ will be $\{x,\emptyset\uplus b_{j_1},\dots,
\emptyset\uplus b_{j_{l-1}} \}$ for some $B=\{b_{j_1},\dots,
b_{j_{l-1}}\}\subseteq V(G)$.
Let $A:=\CN(B)=\{a_{i_1},\dots,a_{i_k}\}$ and
$\tilde\sigma$ be the $\Z_2$-pair of $\sigma$ with
vertex set $\{y,b_{j_1}\uplus\emptyset,\dots,
b_{j_{l-1}}\uplus\emptyset \}$.
We are ready to define $X_{n+1}$:
\[
X_{n+1}:=X_n \setminus \left\{ \tau\in X_n : \sigma\in\tau \mbox{ or }
\tilde\sigma\in\tau\right\}
\]
%
\psfrag{x}{$x$}
\psfrag{B}{$B$}
\psfrag{A}{$A$}
\begin{window}[0,r,{\setlength{\fboxrule}{0pt}\setlength{\fboxsep}{8pt}
\fbox{\includegraphics[width=65mm]{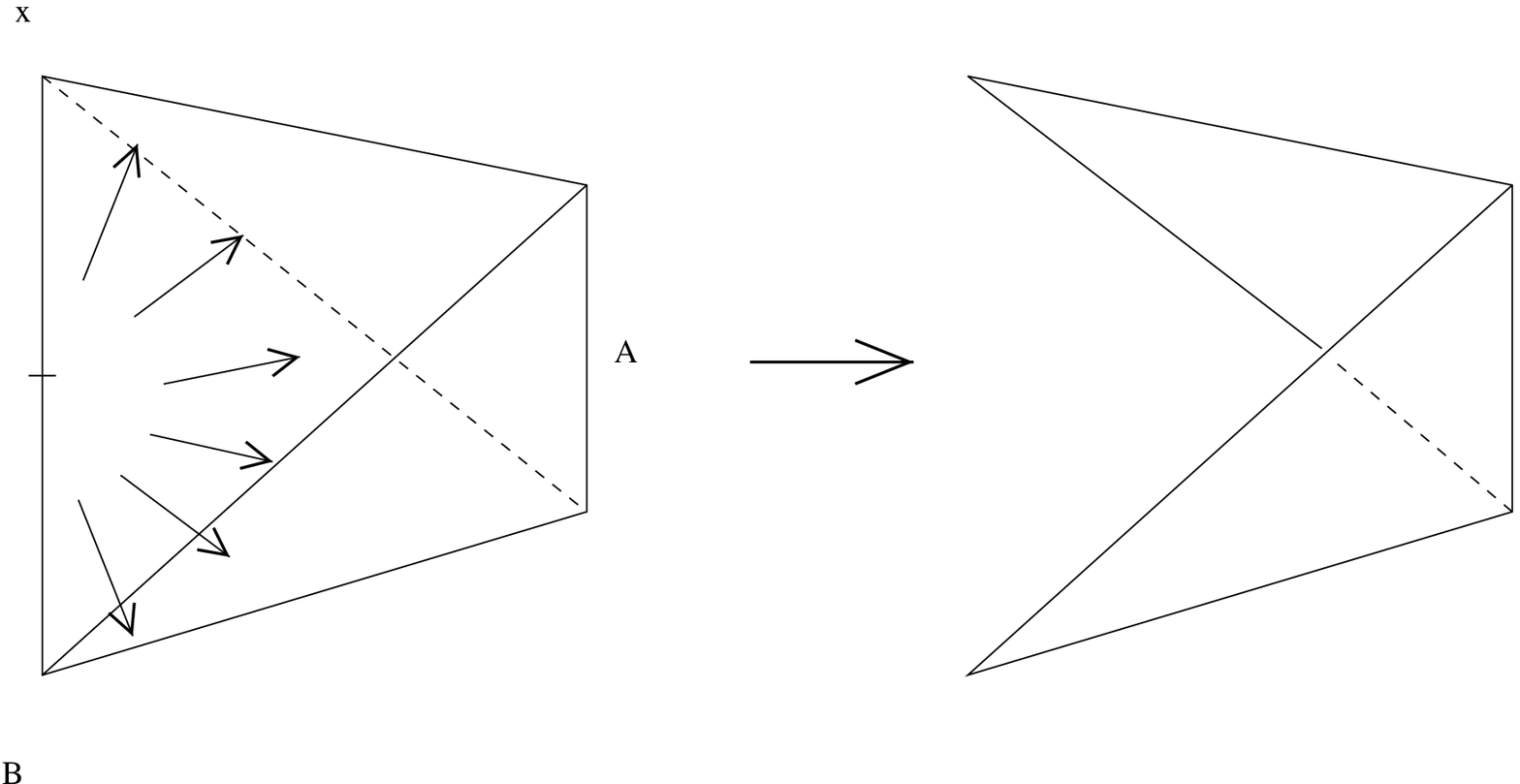} }},{}\label{def}]
We have to only show that $X_{n+1}$ is the deformation retract of $X_n$.
We know the local structure of our complex $X_n$ around $\sigma$. Let
assume that it is a face of a bigger simplex $\Delta$ with vertex set
$\{x,\emptyset\uplus b_{j_1},\dots,\emptyset\uplus b_{j_{l-1}},c\}$.
$c$ can not be the other apex. If $c$ were from the second shore, then we
would choose $\Delta$ instead of $\sigma$ to define $X_{n+1}$. So $c$ can be
only from the first shore and then $c\in A$. This means that $\sigma$
is on the boundary of $X_n$; it is on the boundary of
the simplex $s$ with vertex set $\{x,\emptyset\uplus b_{j_1},\dots,
\emptyset\uplus b_{j_{l-1}},a_{i_1}\uplus\emptyset,\dots,
a_{i_{k}}\uplus\emptyset\}$.
Moreover every simplex which has $\sigma$ as face
is on the boundary of $X_n$. So what we delete to get $X_{n+1}$ is on the
boundary (except $s$). The retraction\footnote{ This deformation
retraction of the simplex $\{v_1=a_{i_1}\uplus\emptyset,\dots,
v_k=a_{i_{k}}\uplus\emptyset,w_1=\emptyset\uplus b_{j_1},\dots,w_{l-1}=
\emptyset\uplus b_{j_{l-1}},w_{l}:=x \}$ can be explicitly given by:
\[
h_t\left( \sum t_i v_i + \sum s_j w_j \right)=
\sum\left(\frac{l\cdot t}{k} +t_i \right)v_i +\sum\left( s_j -t \right) w_j,
\]
where $\sum t_i + \sum s_j =1$. It starts with $h_0=id$,
and ends (for a particular point), just when the first
coefficient of $w_j$ become zero. 
This retraction `kills' those simplices, which has as a face
the simplex $\{w_{1},\dots,w_{l}\}$,
and retracts the `interior' points
to the remaining simplices.}
to $X_{n+1}$ can be given as indicated on the picture.
\hfill $\square$
\end{window}
\end{proofx}

\begin{rem}\label{remark1}
In the same way it can be proven that the neighborhood complex $\N(G)$
(as one shore of the box complex) is a deformation retract of
(homotopy equivalent to) the box complex $\B(G)$.
\end{rem}


\section{Neighborhood complex}\label{sec4}

\begin{figwindow}[0,r,{\setlength{\fboxrule}{0pt}\setlength{\fboxsep}{10pt}
\fbox{\includegraphics[width=35mm]{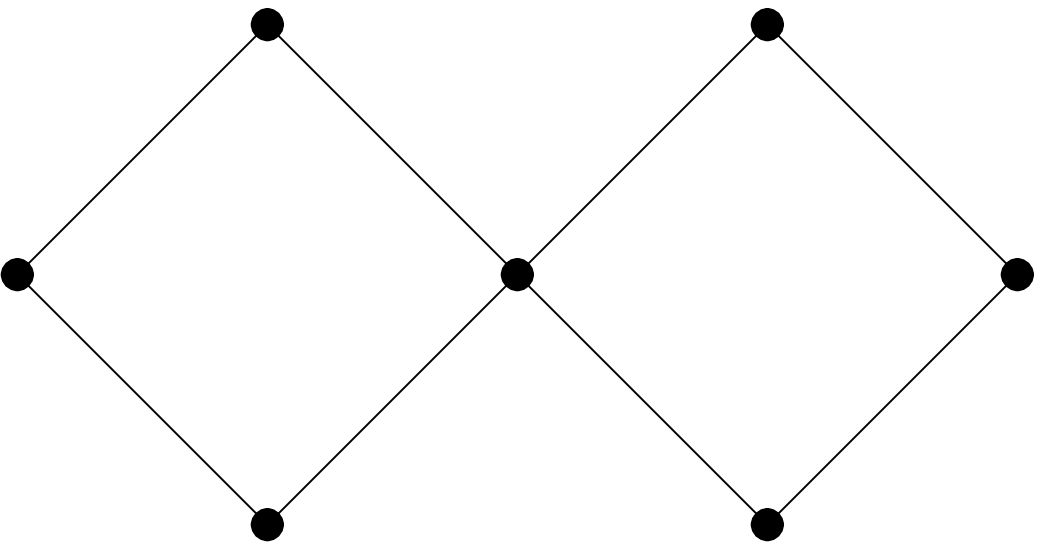} }},{}\label{ex1}]
We consider the following natural question about the neighborhood
complex. Given a simplicial complex $\mathcal{K}$. Is there a graph $G$
such that its neighborhood complex is the given complex,
$\N(G)=\mathcal{K}$?

\indent
For example, if $\mathcal{K}$ is the complex
on Figure \ref{ex1} then the answer is \emph{no}! The reason is that there is a
topological obstruction. The neighborhood complex is homotopy equivalent to
the box complex which is a free $\Z_2$-simplicial complex so it has clearly
even Euler characteristic. But $\mbox{\raisebox{0.8mm}{$\chi$}}(\mathcal{K})=-1$ is odd.
\end{figwindow}
\begin{figwindow}[0,r,{\setlength{\fboxrule}{0pt}\setlength{\fboxsep}{8pt}
\fbox{\includegraphics[width=15mm]{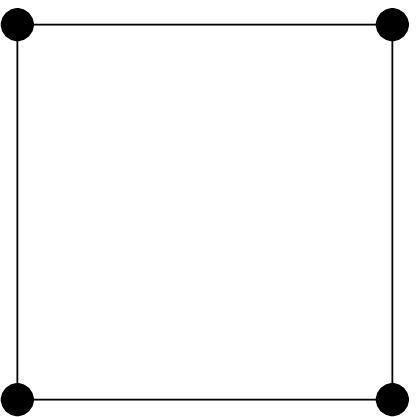} }},{}\label{ex2}]
Another example if $\mathcal{K}$ is the complex
of Figure \ref{ex2}. Now the answer is no again,
but there is no topological reason. With the usual antipodal map $\mathcal{K}$ 
become a free $\Z_2$-simplicial complex. On the other hand the graph $G$
with $\N(G)=\mathcal{K}$ should have 4 vertices, and
by brute force one can
check that $\mathcal{K}$ is not a neighborhood complex.

\indent
Unfortunately we can not answer this question, but we will show that
up to homotopy everything is possible.
\end{figwindow}

\vspace{3mm}
\begin{thm}\label{main}
Given a free $\Z_2$-simplicial complex $\left(\mathcal{K},\nu \right)$,
there is a graph $G$ such that its neighborhood complex is homotopy
equivalent to the given complex,
$\N(G) \simeq \mathcal{K}$.
\end{thm}

In order to prove it we will use the following
construction of a graph from a $\Z_2$-simplicial complex.

\begin{figure}[!h]
\psfrag{G_K}{$G_{\mathcal{K}}$}
\psfrag{K}{$\mathcal{K}$}
\centerline{\includegraphics{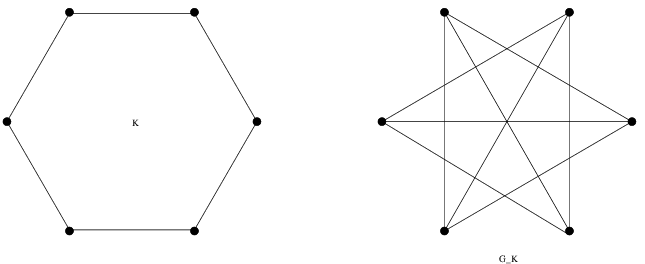}}
\caption{Example for the construction.}\label{ex3}
\end{figure}

\begin{constr}[$\mathcal{K}\to G_{\mathcal{K}}$]\label{con1}
Let $\mathcal{K}$ be a $\Z_2$-simplicial complex. The vertices of
$G_{\mathcal{K}}$ are the vertices of $\mathcal{K}$, and each vertex
is connected to its $\Z_2$-pair and the neighbors\footnote{ in the
1-skeleton of $\mathcal{K}$} of the $\Z_2$-pair.
Thus if $x,y\in V(G_{\mathcal{K}})=V(\mathcal{K})$
then there is an edge between them if and only if $\nu(x)=y$ or
$\{x,\nu(y)\}\in {\mathcal{K}}$ (or $\{y,\nu(x)\}\in {\mathcal{K}}$).
An example is in Figure \ref{ex3}.
\end{constr}

We need the \emph{nerve theorem} as well.

\begin{def2}[Nerve]
Let $\mathcal{F}$ be a set-system. The \emph{nerve} $\mathcal{N}(\mathcal{F})$
of $\mathcal{F}$ is defined as the simplicial complex whose vertices are the
sets in $\mathcal{F}$, and $\left\{ X_1,\dots,X_r \right\}\in
\mathcal{N}(\mathcal{F})$ if and only if $X_1,\dots,X_r \in \mathcal{F}$
and $X_1\cap X_2\cap\cdots\cap X_r\not=\emptyset$.
\end{def2}

\begin{thm}[Nerve theorem \cite{B}]\label{nerve}
Let $\mathcal{K}$ be a simplicial complex and $\mathcal{K}_i$ $(i\in I)$
a family of subcomplexes such that  $\mathcal{K}=\bigcup_{i\in I}\mathcal{K}_i$.
Assume that every nonempty finite intersection $\mathcal{K}_{i_1}\cap
\cdots\cap\mathcal{K}_{i_r}$ is contractible. Then $\mathcal{K}$ and the nerve
$\mathcal{N}(\bigcup\mathcal{K}_i)$ are homotopy equivalent.
\end{thm}


\smallskip
\begin{proof1}
For technical reason we need the first barycentric subdivision
$\sd(\mathcal{K})$ of $\mathcal{K}$. The free simplicial $\Z_2$-action on
$\sd(\mathcal{K})$ will be denoted by $\nu$ as well.

We use Construction \ref{con1} with $\sd(\mathcal{K})$ to obtain
$G_{\sd(\mathcal{K})}$. Because of the barycentric subdivision the vertices
of $G_{\sd(\mathcal{K})}$ denoted by subsets of $V(\mathcal{K})$.
If $A,B\in V(G_{\sd(\mathcal{K})})$
then there is an edge between them if and only if
$\nu(A){=}B$ or $\nu(A){\subset} B$ or $\nu(A){\supset} B$.

We denote the vertices of $\mathcal{K}$ by $1,2,\dots,n$.
Let $\mbox{star}_{\sd(\mathcal{K})}(A)$ be the
star\footnote{ The star of $\sigma\in\mathcal{K}$:
$\mbox{star}_{\mathcal{K}}(\sigma)=\left\{\tau\in\mathcal{K}
\colon \tau\cup\sigma\in\mathcal{K}\right\}$}
of the vertex $A$ in $\sd(\mathcal{K})$. The \emph{nerve} of the set system
$\left\{\mbox{star}_{\sd(\mathcal{K})}(A)
\colon A\in V(G_{\sd(\mathcal{K})})\right\}$
is clearly the neighborhood complex of $G_{\sd(\mathcal{K})}$.
(This is even true without any subdivision:
$N(G_{\mathcal{K}})=\mathcal{N}(\mathcal{S})$ where $\mathcal{S}$ is the
set of the vertex stars in $\mathcal{K}$.)


We want to use the \emph{nerve} theorem so we should prove that if
$B\in\mbox{star}_{\sd(\mathcal{K})}(A_1)\cap\cdots\cap
\mbox{star}_{\sd(\mathcal{K})}(A_r)\not=\emptyset$
then this intersection is contractible. We show that this is a cone.
We have two cases:
\vspace*{-2mm}
\begin{enumerate}
\item If $A_i\subset B$ for all $i=1,2,\dots ,r$.
\\
In this case $\cup A_i$ is a vertex of the barycentric subdivision
since it is a subset of $B$, and it is in the intersection as well.
We show that the intersection can be contracted to this point.
We construct this deformation retraction by letting each vertex to
travel towards $\cup A_i$ with uniform speed. The only thing that
we have to check is that whenever $B_1\subset B_2\subset
\dots\subset B_q$ is a simplex in the intersection, then with the
special vertex $X:=\cup A_{i}$ they form a simplex as well.
First observe that there is an edge between $X$ and $B_l$, $l\in\{1,\dots,q\}$.
If $B_l\subset  A_{i}$ for some $i$ then $B_l\subset X$ as well.
Otherwise $X\subset B_l$.
For the simplex  $B_1\subset B_2\subset\dots\subset B_q$
if $X\subset B_1$ or $X\supset B_q$ then they form a simplex with $X$.
Otherwise there is an
index $k$ such that $B_k\subset X\subset B_{k+1}$. This means that
$B_1, B_2,\dots, B_q,X$ form a simplex.

\item If $B \subset A_{i_j}$ for some $j=1,\dots,k$ ($k\geq 1$),
      and $A_i\subset B$ for the rest.
\\
In this case $B\subset\mathop{\cap}\limits_{j=1}^{k} A_{i_j}\not=\emptyset$ \footnote{ $B\supset\mathop{\cup}\limits_{A_i\subset B} A_{i}$ 
would be good as well, but it can be the emptyset.}
is a vertex of the barycentric
subdivision and the intersection as well.
We show that the intersection can be contracted to this point.
We construct this deformation retraction by letting each vertex
to travel towards $\cap A_{i_j}$ with uniform speed.
We have to show that whenever $B_1\subset B_2\subset
\dots\subset B_q$ is a simplex in the intersection, then with the
special vertex $X:=\cap A_{i_j}$ they form a simplex as well.
First observe that there is an edge between $X$ and $B_l$, $l\in\{1,\dots,q\}$.
If $B_l\supset  A_{i_j}$ for some $i_j$ then $B_l\supset X$ as well.
Otherwise $X\supset B_l$.
For the simplex  $B_1\subset B_2\subset\dots\subset B_q$
if $X\subset B_1$ or $X\supset B_q$ then it is true. Otherwise there is an
index $k$ such that $B_k\subset X\subset B_{k+1}$ which means that
$B_1, B_2,\dots, B_q,X$ form a simplex.

\end{enumerate}
This completes the proof.
\end{proof1}

\section{Box complex}\label{sec5}

In this section we prove our main theorem. It is
the $\Z_2$-extension of Theorem \ref{main}. Later it was proven by Rade
T.\ \v{Z}ivaljevi\'c \cite{Z}.

\begin{thm}\label{boxmain}
Given a free $\Z_2$-simplicial complex $\left(\mathcal{K},\nu \right)$,
there is a graph $G$ such that its box complex $\B(G)$ is $\Z_2$-homotopy
equivalent to the given complex.
\end{thm}

First we need the $\Z_2$-carrier lemma.

\begin{def2}[carrier]
Let $(\mathcal{K},\nu)$ be a $\Z_2$-simplicial complex and $(T,\mu)$
a $\Z_2$-space. A function $C$ taking faces $\sigma$ of $\mathcal{K}$
to subspaces $C(\sigma)$ of $T$, satisfying $C(\nu(\sigma)) = \mu(C(\sigma))$,
is a $\Z_2$-\emph{carrier} if $C(\sigma)\subseteq C(\tau)$ for all
$\sigma\subseteq\tau$.
\end{def2}
\begin{lemma}[$\Z_2$-carrier lemma]\label{carrier}
Assume that for a $\Z_2$-carrier $C$ for any $\sigma\in\mathcal{K}$\ 
$C(\sigma)$ is contractible.
Then any two $\Z_2$-maps $f,g\colon\mathcal{K} \to T$ that
are both carried by $C$ are $\Z_2$-homotopic.
\end{lemma}
\begin{proof}
The proof is straightforward from the definitions.
For details see the proof of Theorem II.9.2 in \cite{LW69}.
\end{proof}

\vspace{4mm}
\begin{proof2}
We will use the same notations as in the proof of Theorem \ref{main}.
Similarly we obtain $G_{\sd(\mathcal{K})}$ by using Construction
\ref{con1}.\ with $\sd(\mathcal{K})$.
We need to show that the box complex $\B(G_{\sd(\mathcal{K})})$ and
$(\mathcal{K},\nu)$ are $\Z_2$-homotopy equivalent.
In order to prove it we will define $\Z_2$-maps
$f\colon\sd(\B(G_{\sd(\mathcal{K})}))\to\sd(\mathcal{K})$ and
$g\colon\sd(\mathcal{K})\to\B(G_{\sd(\mathcal{K})})$. To complete the
proof we will show that $f$ (and $g$) is a $\Z_2$-homotopy equivalence.

The definition of $g$: This is an embedding.
We map a vertex $A\in\sd(\mathcal{K})$ to $A\uplus\emptyset\in
\B(G_{\sd(\mathcal{K})})$ and of course it's $\Z_2$-pair
$\nu(A)\in\sd(\mathcal{K})$ to $\emptyset\uplus A\in
\B(G_{\sd(\mathcal{K})})$. Here we had to choose! If we pick
$\nu(A)$ first than we mapped $\nu(A)$ to $\nu(A)\uplus\emptyset$ and
$A$ to $\emptyset\uplus\nu(A)$. So we have 2 choices for any $\Z_2$-pair
$A,\nu(A)$. This defines a $\Z_2$-map $g$ on the vertex level. We have
to show that $g$ is simplicial. Let $A_1\subset\dots\subset A_l$ be a
simplex $\sigma$ in $\sd(\mathcal{K})$. Since $A_1\uplus\emptyset,\dots,
A_l\uplus\emptyset$, $\emptyset\uplus\nu(A_1),\dots,\emptyset\uplus\nu(A_l)$
form a simplex in $\B(G_{\sd(\mathcal{K})})$ the image of $\sigma$ is a
simplex. (In $G_{\sd(\mathcal{K})}$ $A_i$ is connected to $\nu(A_i)$ and
since $A_i\subset A_j$ or $A_i\supset A_j$ it is connected to $\nu(A_j)$
as well. So $G_{\sd(\mathcal{K})}[\{A_1,\dots,A_l\};\{\nu(A_1),\dots,
\nu(A_l)\}]$ is complete bipartite.)

The definition of $f$: Let $A_1\uplus\emptyset,\dots,A_l\uplus\emptyset$,
$\emptyset\uplus B_1,\dots,\emptyset\uplus B_k$ be the vertices of a
simplex $\sigma$ in $\B(G_{\sd(\mathcal{K})})$. 
$G_{\sd(\mathcal{K})}[\mathcal{A};\mathcal{B}]$ is complete bipartite where
$\mathcal{A}:=\{A_1,\dots,A_l\}$ and $\mathcal{B}:=\{B_1,\dots,B_k\}$.
This means that $\mathcal{A}\subset\mbox{star}_{\sd(\mathcal{K})}
\nu(B_j)$ for any $j\in\{1,\dots,k\}$ so $\mathcal{A}\subset
\mathop{\cap}\limits_{j=1}^{k}\mbox{star}_{\sd(\mathcal{K})}\nu(B_j)$.
From the proof of Theorem \ref{main}.\ we know that
$\mathop{\cap}\limits_{j=1}^{k}\mbox{star}_{\sd(\mathcal{K})}\nu(B_j)$
is a cone with apex $X$.
Since $\mathcal{A},\nu(\mathcal{B})\subset\mbox{star}_{\sd(\mathcal{K})}X$
we have that $Y:=
\mathop{\cap}\limits_{i=1}^{l}\mbox{star}_{\sd(\mathcal{K})}A_i\bigcap
\mathop{\cap}\limits_{j=1}^{k}\mbox{star}_{\sd(\mathcal{K})}\nu(B_j)
\ne\emptyset$.
From the proof of Theorem \ref{main}.\ we know that $Y$ is a cone.
We denote its apex by $X_{\mathcal{A}}^{\mathcal{B}}$ which can be chosen
to be $\mathop{\cap}\limits_{i=1}^{l}A_i\bigcap
\mathop{\cap}\limits_{j=1}^{k}\nu(B_j)$ if it is not the emptyset.
Now we are able to define $f$.
\[
f(\mathcal{A}\uplus\mathcal{B}):=\left\{
\begin{array}{lc}
\mathop{\cap}\limits_{i=1}^{l} A_i\bigcap
\mathop{\cap}\limits_{j=1}^{k} \nu(B_j)
 & \mbox{ if exist,}\\
X_{\mathcal{A}}^{\mathcal{B}} & \mbox{ otherwise.}
\end{array}
\right.
\]
By the construction it is $\Z_2$ on the vertex level. (We can choose
$X_{\mathcal{B}}^{\mathcal{A}}:=\nu(X_{\mathcal{A}}^{\mathcal{B}})$.)
It is simplicial. An edge with two vertices $\mathcal{A}\uplus\mathcal{B}$
and $\tilde{\mathcal{A}}\uplus\tilde{\mathcal{B}}$ ($\tilde{\mathcal{A}}\subset
\mathcal{A}$, $\tilde{\mathcal{B}}\subset\mathcal{B}$) is mapped to two
vertices $S\subset R$ since $X_{\mathcal{A}}^{\mathcal{B}}$ is in the cone
of $X_{\tilde{\mathcal{A}}}^{\tilde{\mathcal{B}}}$. Now a simplex is mapped
to a chain (since every two vertex is comparable by inclusion).

Next we prove that $f\circ \sd(g)\colon\sd(\sd(\mathcal{K}))\to
\sd(\mathcal{K})$ is $\Z_2$-homotopic to $\id_{\mathcal{K}}$.
We will use the $\Z_2$-carrier lemma. We have to construct 'only' a
contractible $\Z_2$-carrier for $f\circ \sd(g)$ and $\id$.
The image of the vertex $v=\{A_1,\dots,A_l\}$, $A_1\subset\dots\subset A_l$ is
$\sd(g)(v)=\{A_{i_1},\dots,A_{i_s}\}\uplus\{\nu(A_{j_1}),\dots,\nu(A_{j_r})\}$.
And now $f(\sd(g)(v))=A_1\cap\dots\cap A_l=A_1$ in this case!
The image of a simplex with vertex set $\{A_{i_1}\},\{A_{i_1},A_{i_2}\},\dots,
\{A_{i_1},\dots,A_{i_l}\}$ is a face of the simplex $A_1\subset\dots\subset 
A_l$. So for a simplex $\sigma\in\sd(\sd(\mathcal{K}))$ with its maximal
vertex $\{A_1,\dots,A_l\}$ we define $C(\sigma):=\{A_1,\dots,A_l\}\in
\sd(\mathcal{K})$. This $C$ is a contractible $\Z_2$-carrier what we need.
$f\circ \sd(g)$ and $\id_{\mathcal{K}}$ are $\Z_2$-homotopic.

Now we show that $g\circ f\colon\sd(\B(G_{\sd(\mathcal{K})}))\to
\B(G_{\sd(\mathcal{K})})$ is $\Z_2$-homotopic to $\id$. Again we construct
a contractible $\Z_2$-carrier for $g\circ f$ and $\id$.
A vertex $\mathcal{A}\uplus\mathcal{B}$ is mapped to
$X_{\mathcal{A}}^{\mathcal{B}}$ by $f$ and to
$X_{\mathcal{A}}^{\mathcal{B}}\uplus\emptyset$ or
$\emptyset\uplus\nu(X_{\mathcal{A}}^{\mathcal{B}})$ by $g\circ f$.
Let $\mathcal{A}_1\uplus\mathcal{B}_1,\dots,\mathcal{A}_n\uplus\mathcal{B}_n$
the vertex set of a simplex $\sigma$ in $\sd(\B(G_{\sd(\mathcal{K})}))$.
($\mathcal{A}_1\subset\dots\subset\mathcal{A}_n$,
$\mathcal{B}_1\subset\dots\subset\mathcal{B}_n$,
$\mathcal{A}_n:=\{A_1,\dots,A_l\}$ and $\mathcal{B}_n:=\{B_1,\dots,B_k\}$).
We consider the subgraph $H$ of $G_{\sd(\mathcal{K})}$ spanned by
$A_1,\dots,A_l,B_1,\dots,B_k$, their $\Z_2$-image under $\nu$ and
$X_{\mathcal{A}_i}^{\mathcal{B}_i},
\nu(X_{\mathcal{A}_i}^{\mathcal{B}_i})$ for any $i\in\{1,\dots,n\}$.
We will use $H$ (actually $\B(H)$) to define the desired carrier.
First of all $\B(H)$ contains the simplex with vertex set
$A_1\uplus\emptyset,\dots,A_l\uplus\emptyset$,
$\emptyset\uplus B_1,\dots,\emptyset\uplus B_k$ which contains $\sigma$.
Moreover we defined $H$ in such a way that $\B(H)$ contains
$(g\circ f)(\sigma)$ as well.
Observe that $H$ is bipartite.
The neighbors of the vertices $X_{\mathcal{A}_n}^{\mathcal{B}_n}$
and 
$\nu(X_{\mathcal{A}_n}^{\mathcal{B}_n})$ provides a partition of
the vertex set of $H$.
The neighborhood complex $\N(H)$ is the disjoint union of two simplices
corresponding to this partition.
So the box complex $\B(H)\subset \B(G_{\sd(\mathcal{K})})$ contains
two disjoint contractible sets (since it is homotopy equivalent to $\N(H)$).
One of these sets covers $\sigma$ and $(g\circ f)(\sigma)$, so we define our
contractible $\Z_2$-carrier $C(\sigma)$ to be this 'half' of $\B(H)$.
\end{proof2}

\begin{rem}
For any free $\Z_2$-simplicial complex $\left(\mathcal{K},\nu \right)$
there is a graph $G$ such that its Hom complex \cite{BK03}
$\mbox{Hom}(K_2,G)$ is $\Z_2$-homotopy equivalent to the given complex
since the box complex $\B(G)$ is $\Z_2$-homotopy
equivalent to $\mbox{Hom}(K_2,G)$.
(The $\Z_2$-maps
$f\colon \sd(\B(G))\to \sd(\mbox{Hom}(K_2,G))$ defined by
$$A\uplus B\to\left\{\begin{array}{lc}
(A,\CN(A)) & \mbox{ if } B=\emptyset,\\
(\CN(B),B) & \mbox{ if } A=\emptyset,\\
(A,B) & \mbox{ otherwise, }
\end{array}
\right.$$ and $g\colon \sd(\mbox{Hom}(K_2,G))\to\sd(\B(G))$
given by $(A,B)\to A\uplus B$
are $\Z_2$-homotopy equivalences. $f\circ g=\mbox{id}$ and
$g\circ f$ is carried by $\mbox{id}$.)
\end{rem}


\section{The suspension and the index}\label{sec6}

In this section we will construct a $\Z_2$-space $X$ such that
$\ind(X)=\ind(\susp(X))$. This example is based on an earlier construction
by Matou\v{s}ek, \v{Z}ivaljevi\'c and the author \cite[page 100]{M}.
Such examples are probably well known for experts (see e.g.\ \cite{CF}),
but we will give a simple and explicit example.

We proceed in the following way.
Let $h\colon S^3\to S^2$ be the Hopf map\footnote{ Considering $S^3$ as the
unit sphere in $\C^2$ and $S^2=\C\P^1$, the Hopf map $h\colon S^3\to S^2$
defined by $(z_1,z_2)\to [z_1,z_2]\in \C\P^1$ \cite[Example 4.45]{Hat01}.
$h$ is a generator of $\pi_3(S^2)=\Z$.}.

We choose a map $2h \colon S^{3} \to S^2$
and we attach two $4$-cells
(the boundary of the $4$-cell is $S^{3}$)
to $S^2$ via $2h$ and $-2h$.
We denote this $\Z_2$-space by
\[X_{2h}:=S^2\mathop{\cup}\limits_{2h} B^{4} \mathop{\cup}_{-2h} B^{4}. \]
The $\Z_2$-action on $S^n\subset X_{2h}$
is the antipodality and interchanges the two $4$-cells.

Now we compute the $\Z_2$-index of $X_{2h}$ and $\susp(X_{2h})$.
%
It is easy to see that $1\leq\ind(X_{2h})\leq 3$.
A $\Z_2$-map $S^2\subset X_{2h} \stackrel{\Z_2}{\longrightarrow} S^1 $
would contradict to the Borsuk--Ulam Theorem.
Let $B^i$ be the unit ball
in $\R^i$ centered at the origin. We assume that $2h\colon S^{3}\to S^2$
maps the unit sphere, the boundary of the unit ball, into the unit sphere.
We define a map $b\colon B^{4}\to B^{3}$ such that it maps the origin of
$\R^{4}$ into the origin of $\R^{3}$ and if $x\in B^{4}$,
$\|x\|\not = 0$ then $b(x) := 2h \left(
\frac{\setlength{\fboxrule}{0pt}\fbox{\scriptsize{$x$}}}{\|x\|}
\right)\cdot \|x\|$. Now we are ready to construct a $\Z_2$-map
$f\colon X_{2h} \stackrel{\Z_2}{\longrightarrow} S^3$.
$f$ maps $S^2\subset X_{2h}$ into the equator of $S^3$.
The remaining two $4$-cells of $X_{2h}$ are mapped to the upper and lower
hemisphere of $S^3$ by $b$ and $-b$.

\medskip
It is slightly more difficult to prove that the index is 3.
We will use the following:
\begin{Def}[\cite{Hat01} Page 427, Section 4.B]
Let $f \colon S^{2n{-}1} \to S^n$, $(n\geq 2)$, and let
$C_f=S^n\mathop{\cup}\limits_f B^{2n}$ (we attach a $2n$-cell to $S^n$ via f).
The Hopf invariant of $f$ (denoted by $\H(f)$) can be defined such that
$\alpha\cup\alpha=\H(f)\cdot\beta$, where $\alpha\in H^n(C_f)=\Z$ and
$\beta\in H^{2n}(C_f)=\Z$  are the generators of the corresponding cohomology
groups and $\cup$ is the cup product.
\end{Def}
We will use the following property of the Hopf invariant (see \cite{Hat01}).
\begin{enumerate}
\item[$\bullet$] $\H:\pi_{2n{-}1}(S^n)\to\Z$ is a homomorphism. For $n{=}2$
it is an isomorphism.
\end{enumerate}
\begin{thm}[\cite{HW60} Theorem 9.5.9]\label{Hopf}
Let $f \colon S^{2n{-}1} \to S^n$ and $g \colon S^n\to S^n$ 
be continuous maps.
Then:
$\H(g \circ f)= \deg (g)^2 \cdot \H(f)$.
\end{thm}

\begin{thm}[\cite{Hat01} Proposition 2B.6]\label{odd}
Every $\Z_2$-map $f \colon S^n {\stackrel{\Z_2}{\longrightarrow}} S^n$
has odd degree.
\end{thm}

\begin{lemma}\label{ind}
$\ind(X_{2h})=3$.
\end{lemma}
\begin{proof}
By contradiction assume that \ $\ind(X_{2h})\leq 2$ which means that there
is a $\Z_2$-map $f\colon X_{2h} \stackrel{\Z_2}{\longrightarrow} S^{2}$.
We restrict this map to $S^2\subset X_{2h}$ obtaining
$g\colon S^2\to S^2$. We claim that $g\circ 2h\colon S^{3}\to S^2$
is null-homotopic.
In $X_{2h}$ we attached a $4$-cell to $S^2$ via
$2h$. This gives us a map $i\colon B^{4}\to X_{2h}$ and
$f\circ i\colon B^{4}\to S^2$. The restriction of $f\circ i$ into
$S^{3}=\partial B^{4}$ is $g\circ 2h$.
So the map $g\circ 2h$ extends into $B^{4}$ which proves
that $g\circ 2h$ is null-homotopic.

On the other hand
Theorem \ref{odd} tells us that $\deg(g)$ is odd.
(We need now only that it is non-zero.)
Using Theorem \ref{Hopf} we
have that $\H(g \circ 2h)= \deg (g)^2 \cdot \H(2h)$.
Since $\deg (g)\not=0$ and $\H(2h)=2$ we have that $\H(g \circ 2h)\not=0$.
This means that $g\circ 2h$ is not null-homotopic, contradiction.
\end{proof}

\begin{lemma}\label{suspind}
$\ind(\susp(X_{2h}))=3$.
\end{lemma}
\begin{proof}
$\susp(X_{2h})$ can be obtained similarly as $X_{2h}$:
we attach two $5$-cells
(the boundary of the $5$-cell is $S^{4}$)
to $S^3$ via $\susp(2h)$ and $-\susp(2h)$.
The Freudenthal Theorem (\cite{Hat01} Corollary 4.24.) tells us that
$\susp\colon\pi_3(S^2)\to\pi_{4}(S^{3})$, which is actually $\Z\to\Z_2$,
is surjective. So $\susp(2h)$ is null-homotopic which means
that $\susp(X_{2h})$ is $\Z_2$-homotopy equivalent to $S^3$ so its
index is 3.
\end{proof}

\bigskip

The generalization of this construction provides infinitely
many examples of $\ind(X)=\ind(\susp(X))$.


Using a simplicial model for $2h\colon S^3_{12}\to S^2_4$
\cite{MS00},\cite{M02} one can obtain a simplicial
complex model for $X_{2h}$ as well.

\section{The topological lower bound can be arbitrarily bad}\label{sec7}

It is well known (see \cite{W}) that the topological lower bound for
the chromatic number 
can be arbitrarily bad. But now we are able
to give purely topological examples.

\begin{def2}
For a graph $G$ let $G^+$ be the graph obtained from $G$ by adding
an extra vertex $w$ and connecting it by edges to all the vertices of $G$,
i.e., $V(G^+)=V(G)\cup \{w\}$ and $E(G^+)=E(G)\cup\{(v,w)\colon v\in V(G)\}$.
\end{def2}

\begin{lemma}
$\B(G^+)$ is $\Z_2$-homotopy equivalent to $\susp(\B(G))$.
\end{lemma}
\begin{proof}
$\susp(\B(G))$ is a subcomplex of $\B(G^+)$. The difference is only
two big simplices (and some of their faces) $V(G)\uplus w$
and $w\uplus V(G)$. We will get rid of the extra simplices one by one
using deformation retraction. We will work with one shore,
on the other shore we have to do the $\Z_2$-pair of each step.

We will define (by induction) sequences of simplicial complexes such that
\[
\B(G^+)=:X_0\supset X_1\supset\dots \supset X_N=\susp(\B(G))
\]
and $X_{i+1}$ is a $\Z_2$-deformation retraction of $X_i$.

Let assume that we already defined $X_n$. We choose $A\subseteq V(G)$
such that $A\uplus w\in X_n$, and there is no $A\subset B\subseteq V(G)$
such that $B\uplus w\in X_n$.
We define $X_{n+1}$:
\[
X_{n+1}:=X_n \setminus \left\{ A\uplus w,w\uplus A,
         A\uplus\emptyset,\emptyset\uplus A \right\}
\]
By the definition of $X_{n+1}$ it is clearly a $\Z_2$-deformation
retract of $X_n$ since $A\uplus \emptyset$ is on the boundary of $X_n$.
(Map the barycenter of $A\uplus \emptyset$ to $\emptyset\uplus w$.)
\end{proof}

\bigskip

Now we are ready to construct a graph such that $\chr G\geq\ind(\B(G))+2+k$.
First we need a $\Z_2$-space (actually a simplicial complex) $X$
such $\ind(X)=\ind(\susp^{k}(X))$.
Now let $G:=G_{sd(X)}$. For $G$ we have that
$\chr G\geq\ind(\B(G))+2=\ind(X)+2$.
We claim that $G^{+k}$ is good for us.
Clearly $\chr G+k=\chr{G^{+k}}$ and $\ind(\B(G^{+k}))=
\ind(\susp^{k}(\B(G)))=\ind(\susp^{k}(X))=\ind(X)$.
So $\chr{G^{+k}}\geq\ind(\B(G^{+k}))+2+k$.

\section{Acknowledgments}
The author would like to thank Ji\v{r}\'\i\ Matou\v{s}ek and
G\"unter M.\ Ziegler for bringing the problems studied in this paper
to my attention and for helpful discussions.



\end{document}